\documentclass[12pt,a4paper]{amsart}

\usepackage{hyperref}

\usepackage{amsrefs}
\usepackage{amssymb}
\usepackage[all]{xy}

\let\mod=\undefined
\DeclareMathOperator{\GL}{GL} %
\DeclareMathOperator{\id}{id} %
\DeclareMathOperator{\op}{op} %
\DeclareMathOperator{\pd}{pd} %
\DeclareMathOperator{\rk}{rk} %
\DeclareMathOperator{\tr}{tr} %
\DeclareMathOperator{\End}{End} %
\DeclareMathOperator{\Ext}{Ext} %
\DeclareMathOperator{\Hom}{Hom} %
\DeclareMathOperator{\mod}{mod} %
\DeclareMathOperator{\rep}{rep} %
\DeclareMathOperator{\supp}{supp} %
\DeclareMathOperator{\gldim}{gldim} %
\DeclareMathOperator{\bdim}{\mathbf{dim}} %

\newcommand{\bd}{\mathbf{d}}

\newcommand{\bh}{\mathbf{h}}

\newcommand{\bbB}{\mathbb{B}}
\newcommand{\bbM}{\mathbb{M}}
\newcommand{\bbN}{\mathbb{N}}
\newcommand{\bbZ}{\mathbb{Z}}

\newcommand{\calC}{\mathcal{C}}
\newcommand{\calE}{\mathcal{E}}
\newcommand{\calF}{\mathcal{F}}
\newcommand{\calH}{\mathcal{H}}
\newcommand{\calO}{\mathcal{O}}
\newcommand{\calP}{\mathcal{P}}
\newcommand{\calQ}{\mathcal{Q}}
\newcommand{\calR}{\mathcal{R}}
\newcommand{\calT}{\mathcal{T}}
\newcommand{\calU}{\mathcal{U}}
\newcommand{\calV}{\mathcal{V}}
\newcommand{\calX}{\mathcal{X}}
\newcommand{\calY}{\mathcal{Y}}

\newcommand{\frakR}{\mathfrak{R}}

\newtheorem{theorem}{Theorem}
\newtheorem{corollary}[theorem]{Corollary}

\newtheorem{proposition}[subsection]{Proposition}

\newcommand{\ol}{\overline}

\title[Irreducible components of module varieties]%
{On regularity in codimension one \\ of irreducible components
\\ of module varieties}

\author{Grzegorz Bobi{\'n}ski}

\address{Faculty of Mathematics and Computer Science \\ Nicolaus
Copernicus University \\ ul.~Chopina 12/18 \\ 87-100 Toru\'n \\
Poland}

\email{gregbob@mat.uni.torun.pl}

\date{}

\keywords{quasi-tilted algebra, module variety, indecomposable
module, irreducible component, regularity in codimension one}

\subjclass[2000]{16G20, 14B05, 14L30}

\begin{document}


\begin{abstract}
Let $\Lambda$ be a tame quasi-tilted algebra and $\bd$ the
dimension vector of an indecomposable $\Lambda$-module. In the
paper we prove that each irreducible component of the variety of
$\Lambda$-modules of dimension vector $\bd$ is regular in
codimension one.
\end{abstract}

\maketitle


Throughout the paper $k$ denotes a fixed algebraically closed
field. Unless otherwise stated all considered modules are finite
dimensional ones. By $\bbZ$, $\bbN$, and $\bbN_+$, we denote the
sets of integers, nonnegative integers, and positive integers,
respectively. For $i, j \in \bbZ$ we denote by $[i, j]$ the set of
all $l \in \bbZ$ such that $i \leq l \leq j$.

\section*{Introduction and the main result}

The class of finite dimensional algebras may be divided into two
disjoint subclasses~\cite{Drozd1980} (see
also~\cite{Crawley-Boevey1988}). The first class consists of the
tame algebras for which indecomposable modules occur in each
dimension in a finite number of one-parameter
families~\eqref{subsect_tame}. The second class is formed by the
wild algebras whose representation theory is as complicated as the
study of finite dimensional vector spaces together with two (not
necessarily commuting) endomorphisms. The solution of the latter
problem would imply the classification of modules over arbitrary
finite dimensional algebra, hence one may hope to solve
classification problems only for tame algebras.

Given a finite dimensional algebra $\Lambda$ and a nonnegative
element $\bd$ of the Grothendieck group of the category of
$\Lambda$-modules (such elements of the Grothendieck group are
called dimension vectors), it is an interesting problem to study
the variety $\mod_\Lambda (\bd)$ of $\Lambda$-modules of dimension
vector $\bd$~\eqref{subsect_modulevariety} possessing a natural
action of a corresponding product $\GL (\bd)$ of general linear
groups~\eqref{subsect_action} (see~\cites{Bobinski2002,
Bobinski2008a, BobinskiSkowronski1999a, BobinskiSkowronski1999b,
BobinskiSkowronski2002, Bongartz1991, Crawley-BoeveySchroer2002,
Geiss1996, GeissSchroer2003, GeissSchroer2005, Richmond2001,
Schroer2004} for some results in this direction). In particular,
one may ask questions about irreducible components of
$\mod_\Lambda (\bd)$ and their properties. This problem is
nontrivial only for algebras which are not hereditary, since
$\mod_\Lambda (\bd)$ is an affine space provided $\gldim \Lambda
\leq 1$. A natural class of algebras obtained by weakening
homological constrains put on hereditary algebras is that of the
quasi-tilted algebras~\eqref{subsect_quasitilted}. It seems to be
natural to study how much varieties of modules over quasi-tilted
algebras differ from varieties of modules over hereditary
algebras.

The categories of modules over quasi-tilted algebras have been
intensively studied (see for example~\cites{HappelReiten1999,
HappelReitenSmalo1996, HappelRingel1982, LenzingdelaPena1997,
Ringel1984}). In particular, the categories of modules over tame
quasi-tilted algebras have been described
(see~\cite{Skowronski1998}, also~\cites{Kerner1989,
LenzingSkowronski1996, delaPena1993}). Since some geometric
properties of points of varieties of modules translate into
properties of corresponding modules, it is appealing to use
information about categories of modules in studying varieties of
modules.

This strategy has been successfully applied
in~\cites{Bobinski2002, BobinskiSkowronski1999a,
BobinskiSkowronski1999b, BobinskiSkowronski2002}. In particular,
the main results of~\cites{BobinskiSkowronski1999a,
BobinskiSkowronski1999b} imply that if $\bd$ is the dimension
vector of an indecomposable module over a tame quasi-tilted
algebra $\Lambda$, then $\mod_\Lambda (\bd)$ has at most two
irreducible components. Moreover, if $\mod_\Lambda (\bd)$ is
irreducible, then it is normal, hence its unique irreducible
component is regular in codimension one~\eqref{subsect_codimone}.
However, no further information about properties of the
irreducible components of $\mod_\Lambda (\bd)$ in the remaining
cases was obtained at that time. A first step in order to fill
this gap was made in~\cite{BobinskiZwara2006} by proving that if
$X$ is an indecomposable directing
module~\eqref{subsect_directing} over a tame algebra $\Lambda$
(which may be assumed without loss of generality to be
quasi-tilted) of dimension vector $\bd$, then the closure of the
$\GL (\bd)$-orbit of $X$ (which is an irreducible component of
$\mod_\Lambda (\bd)$) is normal, hence again regular in
codimension one. The main result of the paper extends the latter
property to all irreducible components of varieties of modules
over quasi-tilted algebras for dimension vectors of indecomposable
modules.

\begin{theorem}
Let $\Lambda$ be a quasi-tilted algebra and $\bd$ a dimension
vector of an indecomposable $\Lambda$-module. If $\calV$ is an
irreducible component of $\mod_\Lambda (\bd)$, then $\calV$ is
regular in codimension one.
\end{theorem}

We remark that the dimension vectors of the indecomposable modules
over the tame quasi-tilted algebras can be characterized
combinatorially~\eqref{subsect_quasitilted}.

An important direction of research is study of properties of
closures of orbits in varieties of modules (see for
example~\cites{BenderBongartz2003, Bongartz1994,
BobinskiZwara2002, SkowronskiZwara2003, Zwara2005b, Zwara2005,
Zwara2006, Zwara2007}). If $X$ is a module of dimension vector
$\bd$ over an algebra $\Lambda$ with $\Ext_\Lambda^1 (X, X) = 0$,
then one easily gets that the closure of the $\GL (\bd)$-orbit of
$X$ is an irreducible component of $\mod_\Lambda
(\bd)$~\cite{Kraft1984}*{II.2.7, Folgerung}. Thus we have the
following immediate consequence of the main result of the paper.

\begin{corollary}
Let $X$ be an indecomposable module over a tame quasi-tilted
algebra $\Lambda$ of dimension vector $\bd$. If $\Ext_\Lambda^1
(X, X) = 0$, then the closure of the $\GL (\bd)$-orbit of $X$ is
regular in codimension one.
\end{corollary}

It is an open problem if the closure of the orbit of $X$ is
regular in codimension one for an arbitrary indecomposable module
$X$ over a tame quasi-tilted algebra. In~\cite{Zwara2003} Zwara
presented an example showing that in general closures of orbits of
(arbitrary) modules over tame quasi-tilted (even hereditary)
algebras do not have to be regular in codimension one.

The paper is organized as follows. In Section~\ref{sect_quivers}
we present necessary information about quivers and quasi-tilted
algebras, while in Section~\ref{sect_varieties} we collect useful
facts about varieties of modules. Next, in
Section~\ref{sect_bisection} we describe technical conditions
under which a distinguished component of a module variety is
regular in codimension one. This condition is used among others in
Section~\ref{sect_proof} in order to prove the main result.

For a basic background on the representation theory of algebras
(in particular, on tilting theory) we refer
to~\cite{AssemSimsonSkowronski2006}. Basic algebraic geometry used
in the article can be found in~\cite{Kunz1985}.

The article was written while the author was staying at the
University of Bielefeld as the Alexander von Humboldt Foundation
fellow.

\makeatletter
\def\@secnumfont{\mdseries} 
\makeatother

\section{Preliminaries on quivers and quasi-tilted algebras}
\label{sect_quivers}

\subsection{} \label{subsect_tame}
A finite dimensional algebra $\Lambda$ is called tame if for every
$d \in \bbN_+$ there exist $\Lambda$-$k [T]$-bimodules $M_1$,
\ldots, $M_n$, which are free of rank $d$ as $k [T]$-modules (in
particular, in contrary to our usual assumption, they are not
finite dimensional) such that for each indecomposable
$\Lambda$-module $X$ there exist $i \in [1, n]$ and $\lambda \in
k$ such that $X \simeq M_i \otimes_{k [T]} (k [T] / (T -
\lambda))$.

\subsection{}
By a quiver $\Delta$ we mean a finite set $\Delta_0$ of vertices
and a finite set $\Delta_1$ of arrows together with maps $s, t :
\Delta_1 \to \Delta_0$ which assign to $\alpha \in \Delta_1$ the
starting vertex $s_\alpha$ and the terminating vertex $t_\alpha$.
By a path of length $n \in \bbN_+$ we mean a sequence $\sigma =
\alpha_1 \cdots \alpha_n$ with $\alpha_1, \ldots, \alpha_n \in
\Delta_1$ such that $s_{\alpha_i} = t_{\alpha_{i + 1}}$ for each
$i \in [1, n - 1]$. In the above situation we put $s_\sigma =
s_{\alpha_n}$ and $t_\sigma = t_{\alpha_1}$. Additionally, for
each $x \in \Delta_0$ we introduce the trivial path of length $0$
starting and terminating at $x$ and denoted by $x$. A path
$\sigma$ of positive length such that $s_\sigma = t_\sigma$ is
called an oriented cycle. A subquiver $\Delta'$ of $\Delta$ (i.e.\
a pair $(\Delta_0', \Delta_1')$ such that $\Delta_0' \subseteq
\Delta_0$, $\Delta_1' \subseteq \Delta_1$, and $s_\alpha, t_\alpha
\in \Delta_0'$ for each $\alpha \in \Delta_1'$) is called convex
if for every path $\alpha_1 \cdots \alpha_n$ in $\Delta$ with
$\alpha_1, \ldots, \alpha_n \in \Delta_1$ and $s_{\alpha_n},
t_{\alpha_1} \in \Delta_0'$, $s_{\alpha_i} \in \Delta_0'$ for each
$i \in [1, n - 1]$.

\subsection{}
For a quiver $\Delta$ we denote by $k \Delta$ the path algebra of
$\Delta$ defined as follows. The elements of $k \Delta$ are the
formal linear combinations of paths in $\Delta$ and for two paths
$\sigma_1$ and $\sigma_2$ the product of $\sigma_1$ and $\sigma_2$
is either the composition $\sigma_1 \sigma_2$ of paths, if
$s_{\sigma_1} = t_{\sigma_2}$, or $0$, otherwise. Fix $x, y \in
\Delta_0$ and let $\rho = \lambda_1 \sigma_1 + \cdots + \lambda_n
\sigma_n$ for $n \in \bbN_+$, $\lambda_1, \ldots, \lambda_n \in
k$, and paths $\sigma_1$, \ldots, $\sigma_n$. If $s_{\sigma_i} =
x$ and $t_{\sigma_i} = y$ for each $i \in [1, n]$, then we write
$s_{\rho} = x$ and $t_{\rho} = y$. If additionally, the length of
$\sigma_i$ is at least (more than) $1$ for each $i \in [1, n]$,
then $\rho$ is called a relation (an admissible relation,
respectively). A set $\frakR$ of relations is called minimal if
$\rho$ does not belong to the ideal $\langle \frakR \setminus \{
\rho \} \rangle$ generated by $\frakR \setminus \{ \rho \}$ for
each $\rho \in \frakR$. If $\frakR$ is a minimal set of admissible
relations such that there exists $n \in \bbN_+$ with the property
$\sigma \in \langle \frakR \rangle$ for each path $\sigma$ of
length at least $n$, then the pair $(\Delta, \frakR)$ is called a
bound quiver. If $(\Delta, \frakR)$ is a bound quiver, then $k
\Delta / \langle \frakR \rangle$ is called the path algebra of
$(\Delta, \frakR)$.

Gabriel proved (see for
example~\cite{AssemSimsonSkowronski2006}*{Corollaries~I.6.10
and~II.3.7}) that for each finite dimensional algebra $\Lambda$
there exists a bound quiver $(\Delta, \frakR)$ such that the
category $\mod \Lambda$ of $\Lambda$-modules is equivalent to the
category of modules over the path algebra of $(\Delta, \frakR)$.
In addition, $\Delta$ is uniquely (up to isomorphism) determined
by $\Lambda$ and we call it the Gabriel quiver of $\Lambda$.
Consequently, from now on we assume that all considered algebras
are path algebras of bound quivers. In particular, all considered
algebras will be finite dimensional.

\subsection{}
Let $\Delta$ be a quiver. By a representation of $\Delta$ we mean
a collection $M = (M_x, M_\alpha)_{x \in \Delta_0, \alpha \in
\Delta_1}$ of finite dimensional vector spaces $M_x$, $x \in
\Delta_0$, and linear maps $M_\alpha : M_{s_\alpha} \to
M_{t_\alpha}$, $\alpha \in \Delta_1$. If $M$ and $N$ are
representations of $\Delta$, then the morphism space $\Hom_\Delta
(M, N)$ consists of the collections $f = (f_x)_{x \in \Delta_0}$
of linear maps $f_x : M_x \to N_x$, $x \in \Delta_0$, such that
$N_{\alpha} f_{s_\alpha} = f_{t_\alpha} M_\alpha$ for each $\alpha
\in \Delta_1$. If $\sigma = \alpha_1 \cdots \alpha_n$, for $n \in
\bbN_+$ and $\alpha_1, \ldots, \alpha_n \in \Delta_1$, is a path,
then for a representation $M$ of $\Delta$ we put $M_\sigma =
M_{\alpha_1} \cdots M_{\alpha_n}$. Similarly, if $\rho = \lambda_1
\sigma_1 + \cdots + \lambda_n \sigma_n$, for $n \in \bbN_+$,
$\lambda_1, \ldots, \lambda_n \in k$, and paths $\sigma_1$,
\ldots, $\sigma_n$ in $\Delta$, is a relation, then for a
representation $M$ of $\Delta$ we put $M_\rho = \lambda_1
M_{\sigma_1} + \cdots + \lambda_n M_{\sigma_n}$.

\subsection{}
Let $(\Delta, \frakR)$ be a bound quiver. By $\rep (\Delta,
\frakR)$ we denote the full subcategory of the category of
representations of $\Delta$ consisting of the representations $M$
such that $M_{\rho} = 0$ for each $\rho \in \frakR$. If $\Lambda$
is the path algebra of $(\Delta, \frakR)$, then the assignment
which assigns to a $\Lambda$-module $M$ the representation $(M_x,
M_\alpha)_{x \in \Delta_0, \alpha \in \Delta_1}$, where $M_x = x
M$ for $x \in \Delta_0$ and $M_\alpha (m) = \alpha m$ for $\alpha
\in \Delta_1$ and $m \in M_{s_\alpha}$, induces an equivalence
between $\mod \Lambda$ and $\rep (\Delta, \frakR)$ (see for
example~\cite{AssemSimsonSkowronski2006}*{Theorem~III.1.6}). We
will usually treat this equivalence as identification. With $M \in
\rep (\Delta, \frakR)$ (hence also with every $\Lambda$-module) we
associate its dimension vector $\bdim M \in \bbZ^{\Delta_0}$
defined by $(\bdim M)_x = \dim_k M_x$ for $x \in \Delta_0$.
Obviously, $\bdim M \in \bbN^{\Delta_0}$. We call the elements of
$\bbN^{\Delta_0}$ dimension vectors. By the support of a dimension
vector $\bd$ we mean the full subquiver $\supp \bd$ of $\Delta$
with $\{ x \in \Delta_0 \mid d_x \neq 0 \}$ as the set of vertices
(a subquiver $\Delta'$ of $\Delta$ is called full if $\Delta_1'$
consists precisely of all $\alpha \in \Delta_1$ such that
$s_\alpha, t_\alpha \in \Delta_0'$). A dimension vector $\bd$ is
called sincere (connected) if $\supp \bd = \Delta_0$ ($\supp \bd$
is connected, respectively).

\subsection{}
Let $(\Delta, \frakR)$ be a bound quiver. With $(\Delta, \frakR)$
we associate the bilinear form $\langle -, - \rangle_{\Delta,
\frakR} : \bbZ^{\Delta_0} \times \bbZ^{\Delta_0} \to \bbZ$ and the
quadratic form $q_{\Delta, \frakR} : \bbZ^{\Delta_0} \to \bbZ$ in
the following way:
\[
\langle \bd', \bd'' \rangle_{\Delta, \frakR} = \sum_{x \in
\Delta_0} d_x' d_x'' - \sum_{x \in \Delta_1} d_{s_\alpha}'
d_{t_\alpha}'' + \sum_{\rho \in \frakR} d_{s_\rho}' d_{t_\rho}''
\]
for $\bd', \bd'' \in \bbZ^{\Delta_0}$ and $q_{\Delta, \frakR}
(\bd) = \langle \bd, \bd \rangle_{\Delta, \frakR}$ for $\bd \in
\bbZ^{\Delta_0}$. Let $\Lambda$ be the path algebra of $(\Delta,
\frakR)$. If $\Lambda$ is triangular (i.e.\ there are no oriented
cycles in $\Delta$), then $\langle -, - \rangle_{\Delta', \frakR'}
= \langle -, - \rangle_{\Delta, R}$ for each bound quiver
$(\Delta', \frakR')$ such that $k \Delta' / \langle \frakR'
\rangle \simeq \Lambda$~\cite{Bongartz1983}*{Proposition~1.2}.
Thus in this case we may write $\langle -, - \rangle_\Lambda$ and
$q_\Lambda$ instead of $\langle -, - \rangle_{\Delta, \frakR}$ and
$q_{\Delta, \frakR}$, respectively. Additionally, if $\gldim
\Lambda \leq 2$, then
\begin{multline*}
\langle \bdim M, \bdim N \rangle_\Lambda
\\
= \dim_k \Hom_\Lambda (M, N) - \dim_k \Ext_\Lambda^1 (M, N) +
\dim_k \Ext_\Lambda^2 (M, N)
\end{multline*}
for all $\Lambda$-modules $M$ and
$N$~\cite{Bongartz1983}*{Proposition~2.2}.

\subsection{} \label{subsect_quasitilted}
An algebra $\Lambda$ is called quasi-tilted if $\gldim \Lambda
\leq 2$ and either $\pd_\Lambda X \leq 1$ or $\id_\Lambda X \leq
1$ for each indecomposable $\Lambda$-module $X$. Equivalently, an
algebra $\Lambda$ is quasi-tilted if and only if there exists a
tilting object $T$ in a hereditary category $\calH$ such that
$\Lambda \simeq \End_{\calH}
(T)^{\op}$~\cite{HappelReitenSmalo1996}*{Theorem~II.2.3}. The
quasi-tilted algebras are
triangular~\cite{HappelReitenSmalo1996}*{Proposition~III.1.1~(b)},
hence the forms $\langle -, - \rangle_\Lambda$ and $q_\Lambda$ are
defined for a quasi-tilted algebra $\Lambda$. Moreover, a
quasi-tilted algebra $\Lambda$ is tame if and only if $q_\Lambda
(\bd) \geq 0$ for each dimension vector
$\bd$~\cite{Skowronski1998}*{Theorem~A}. Additionally, if $\bd$ is
a dimension vector, then there exists an indecomposable
$\Lambda$-module of dimension vector $\bd$ if and only if $\bd$ is
connected and $q_\Lambda (\bd) \in \{ 0, 1 \}$. Finally, if $\bd$
is a connected dimension vector and $q_\Lambda (\bd) = 1$, then
there exists a unique (up to isomorphism) indecomposable
$\Lambda$-module of dimension vector $\bd$.

\subsection{} \label{subsect_tilted} \label{subsect_directing}
A important class of quasi-tilted algebra is that of the tilted
algebras. An algebra $\Lambda$ is called tilted if there exists a
tilting module $T$ over a hereditary algebra $H$ such that
$\Lambda \simeq \End_H (T)^{\op}$ (see for
example~\cite{AssemSimsonSkowronski2006}*{Chapter~VIII} for more
on tilted algebras). Recall that a module $T$ over an algebra
$\Lambda$ is called tilting if $\pd_\Lambda T \leq 1$,
$\Ext_\Lambda^1 (T, T) = 0$, and $T$ is a direct sum of $n$
pairwise nonisomorphic indecomposable $\Lambda$-modules, where $n$
is the number of the vertices in the Gabriel quiver of $\Lambda$.
For a tilting module $T$ over an algebra $\Lambda$ we define two
full subcategories $\calF (T)$ and $\calT (T)$ of $\mod \Lambda$
by
\begin{gather*}
\calF (T) = \{ M \in \mod \Lambda \mid \Hom_\Lambda (T, M) = 0 \}
\\
\intertext{and} %
\calT (T) = \{ N \in \mod \Lambda \mid \Ext_\Lambda^1 (T, N) = 0
\}.
\end{gather*}
Bakke proved~\cite{Bakke1988}*{Theorem} that an algebra $\Lambda$
is tilted if and only if there exists a directing tilting
$\Lambda$-module $T$. Here, a module $T$ over an algebra $\Lambda$
is called directing if and only if there exists no sequence
\[
X_0 \xrightarrow{f_1} X_1 \to \cdots \to X_{n - 1}
\xrightarrow{f_n} X_n
\]
of nonzero nonisomorphism $f_1$, \ldots, $f_n$ between
indecomposable $\Lambda$-modules $X_0$, \ldots, $X_n$ such that
$X_0$ and $X_n$ are direct summands of $T$ and $X_{i - 1} \simeq
\tau_\Lambda X_{i + 1}$ for some $i \in [1, n - 1]$, where
$\tau_\Lambda$ denotes the Auslander--Reiten translation of $\mod
\Lambda$ (see~\cite{AssemSimsonSkowronski2006}*{Section~IV.2} for
a definition and basic properties of the Auslander--Reiten
translation). The following properties of directing tilting
modules will be useful for us
(see~\citelist{\cite{AssemSimsonSkowronski2006}*{Chapter~VIII}
\cite{Bakke1988}}): if $T$ is a directing tilting module over an
algebra $\Lambda$, then $\mod \Lambda = \calF (T) \vee \calT (T)$,
$\pd_\Lambda M \leq 1$ for each $M \in \calF (T)$, $\id_\Lambda N
\leq 1$ for each $N \in \calT (T)$, and $\Hom_\Lambda (N, M) = 0$
and $\Ext_\Lambda^1 (M, N) = 0$ for all $M \in \calF (T)$ and $N
\in \calT (T)$. In the paper, if $\calX$ and $\calY$ are two full
subcategories of $\mod \Lambda$ for an algebra $\Lambda$, then we
denote by $\calX \vee \calY$ the additive closure of their union.

\subsection{} \label{subsect_concan}
Another important class of quasi-tilted algebras is formed by the
concealed canonical algebras introduced
in~\cite{LenzingMeltzer1996} --- we refer
to~\cite{LenzingMeltzer1996} for the definition of concealed
canonical algebras and proofs of their properties listed below. If
$\Lambda$ is a concealed canonical algebra, then it possesses a
sincere separating tubular family $\calR$ (in the sense
of~\cite{Ringel1984}, see also~\cites{LenzingdelaPena1999,
Skowronski1996}). Moreover, if
\begin{gather*}
\calP = \{ M \in \mod \Lambda \mid \text{$\Hom_\Lambda (R, M) = 0$
for each $R \in \calR$} \}
\\
\intertext{and} %
\calQ = \{ N \in \mod \Lambda \mid \text{$\Hom_\Lambda (N, R) = 0$
for each $R \in \calR$} \},
\end{gather*}
then $\mod \Lambda = \calP \vee \calR \vee \calQ$, $\pd_\Lambda M
\leq 1$ for each $M \in \calP \vee \calR$, $\id_\Lambda N \leq 1$
for each $N \in \calR \vee \calQ$, and $\Hom_\Lambda (N, M) = 0$
and $\Ext_\Lambda^1 (M, N) = 0$ if either $M \in \calP \vee \calR$
and $N \in \calQ$ or $M \in \calP$ and $N \in \calR \vee \calQ$.

\makeatletter
\def\@secnumfont{\mdseries} 
\makeatother

\section{Module varieties} \label{sect_varieties}

\subsection{} \label{subsect_modulevariety}
Let $\Lambda$ be the path algebra of a bound quiver $(\Delta,
\frakR)$ and $\bd \in \bbN^{\Delta_0}$. By the variety
$\mod_\Lambda (\bd)$ of $\Lambda$-modules of dimension vector
$\bd$ we mean the set consisting of all $M \in \rep (\Delta,
\frakR)$ such that $M_x = k^{d_x}$. By forgetting the spaces
$M_x$, $x \in \Delta$, we identify $\mod_\Lambda (\bd)$ with the
Zariski-closed subset of $\prod_{\alpha \in \Delta_1}
\bbM_{d_{t_\alpha} \times d_{s_\alpha}} (k)$. Observe that for
each $\Lambda$-module $N$ of dimension vector $\bd$ there exists
$M \in \mod_\Lambda (\bd)$ such that $M \simeq N$. We will usually
assume that all considered $\Lambda$-modules are points of module
varieties. Observe that the above construction can be performed
for any pair $(\Delta', \frakR')$ consisting of a quiver $\Delta'$
and a set of relations $\frakR'$ (not only for bound quivers). The
varieties obtained in this way are called representation
varieties.

\subsection{} \label{subsect_deg} \label{subsect_action}
Let $\Lambda$ be an algebra, $\Delta$ its Gabriel quiver, and $\bd
\in \bbN^{\Delta_0}$. Then $\GL (\bd) = \prod_{x \in \Delta_0}
\GL_{d_x} (k)$ acts on $\mod_\Lambda (\bd)$ by conjugation: $(g
M)_{\alpha} = g_{t_\alpha} M_\alpha g_{s_\alpha}^{-1}$ for $g \in
\GL (\bd)$, $M \in \mod_\Lambda (\bd)$, and $\alpha \in \Delta_1$.
Observe that $M \simeq N$ for $M, N \in \mod_\Lambda (\bd)$ if and
only if $\calO (M) = \calO (N)$, where for $M \in \mod_\Lambda
(\bd)$ we denote by $\calO (M)$ its $\GL (\bd)$-orbit. One easily
calculates that $\dim \calO (M) = \dim \GL (\bd) - \dim_k
\End_\Lambda (M)$ for each $M \in \mod_\Lambda
(\bd)$~\cite{KraftRiedtmann1986}*{2.2}. If $\calU$ is a $\GL
(\bd)$-invariant subset of $\mod_\Lambda (\bd)$ and $M \in \calU$,
then we say that $\calO (M)$ is maximal in $\calU$ if there is no
$N \in \calU$ such that $\calO (M) \subseteq \ol{\calO (N)}$ and
$N \not \simeq M$.

Fix $M, N \in \mod_\Lambda (\bd)$. The formula for the dimension
of an orbit implies in particular that $\dim_k \End_\Lambda (M) <
\dim_k \End_\Lambda (N)$ provided $\calO (N) \subseteq \ol{\calO
(M)}$ and $N \not \simeq M$. Moreover, if there exists an exact
sequence $0 \to N_1 \to M \to N_2 \to 0$ such that $N_1 \oplus N_2
\simeq N$, then $\calO (N) \subseteq \ol{\calO (M)}$ (see for
example~\cite{Bongartz1996}*{Lemma~1.1}). On the other hand, it
has been proved in~\cite{Zwara2000}*{Corollary~6} that if
$\Lambda$ is tame and quasi-tilted, then $\calO (N) \subseteq
\ol{\calO (M)}$ if and only if there exist exact sequences $0 \to
U_i' \to V_i \to U_i'' \to 0$, $i \in [0, n]$, for $n \in \bbN$,
such that $M \simeq V_0$, $V_i \simeq U_{i - 1}' \oplus U_{i -
1}''$ for each $i \in [1, n]$, and $U_n' \oplus U_n'' \simeq N$.

\subsection{} \label{subsect_codimone}
Let $\calV$ be an affine variety. Recall that $x \in \calV$ is
called a regular point of $\calV$ if the dimension of the tangent
space $T_x \calV$ to $\calV$ at $x$ equals the maximum of the
dimensions of the irreducible components of $\calV$ containing $x$
(see for example~\cite{Kunz1985}*{VI.1}). In particular, if
$\calV$ is irreducible, then $x \in \calV$ is a regular point of
$\calV$ if and only if $\dim_k T_x \calV = \dim \calV$. In
general, if $\calV'$ is an irreducible component of $\calV$ and $x
\in \calV'$ is a regular point of $\calV$, then $x$ does not
belong to an irreducible component of $\calV$ different from
$\calV'$. In particular, in the above situation $x$ is a regular
point of $\calV'$. We say that $\calV$ is regular in codimension
one if the codimension of the complement of the set of regulars
points of $\calV$ is at least $2$. For example, normal varieties
are regular in codimension one~\cite{Eisenbud1995}*{Section~11.2}.
Recall, that $\calV$ is called normal if $\calV$ is irreducible
and its coordinate ring is an integrally closed domain.

\subsection{} \label{subsect_cocyle}
Let $\Lambda$ be the path algebra of a bound quiver $(\Delta,
\frakR)$ and $\bd', \bd'' \in \bbN^{\Delta_0}$. Fix $U \in
\mod_\Lambda (\bd')$ and $V \in \mod_\Lambda (\bd'')$. If $\rho =
\lambda_1 \sigma_1 + \cdots + \lambda_n \sigma_n$, for $n \in
\bbN_+$, $\lambda_1, \ldots, \lambda_n \in k$, and paths
$\sigma_1$, \ldots, $\sigma_n$ in $\Delta$, is a relation, then
for $Z \in \prod_{\alpha \in \Delta_1} \bbM_{d_{t_\alpha}' \times
d_{s_\alpha}''} (k)$ we put
\[
Z_\rho = \sum_{i \in [1, n]} \sum_{j \in [1, m_i]} \lambda_i
U_{\alpha_{i, 1}} \cdots U_{\alpha_{i, j - 1}} Z_{\alpha_{i, j}}
V_{\alpha_{i, j + 1}} \cdots V_{\alpha_{i, m_i}},
\]
provided $\sigma_i = \alpha_{i, 1} \cdots \alpha_{i, m_i}$ with
$\alpha_{i, 1}, \ldots, \alpha_{i, m_i} \in \Delta_1$ and $m_i \in
\bbN_+$ for each $i \in [1, n]$. Let $\bbZ (V, U)$ be the set of
all $Z \in \prod_{\alpha \in \Delta_1} \bbM_{d_{t_\alpha}' \times
d_{s_\alpha}''} (k)$ such that $Z_\rho = 0$ for each $\rho \in
\frakR$. With every $Z \in \bbZ (V, U)$ we associate the exact
sequence $\xi^Z : 0 \to U \xrightarrow{f} W^Z \xrightarrow{g} V
\to 0$ of $\Lambda$-modules, where $W_x^Z = U_x \oplus V_x$ for $x
\in \Delta_0$, $W_\alpha = \left[
\begin{smallmatrix}
U_\alpha & Z_\alpha \\ 0 & V_\alpha
\end{smallmatrix}
\right]$ for $\alpha \in \Delta_1$, and $f_x$ and $g_x$ are the
canonical injection and projection, respectively, for $x \in
\Delta_0$. This assignment induces a surjective linear map $\bbZ
(V, U) \to \Ext_\Lambda^1 (V, U)$. Let $\bbB (V, U)$ be the kernel
of this map. Then
\begin{gather*}
\dim_k \bbB (V, U) =  - \dim_k \Hom_\Lambda (V, U) + \sum_{x \in
\Delta_0} d_x' d_x''
\\
\intertext{and, consequently,} %
\dim_k \bbZ (V, U) = \dim_k \Ext_\Lambda^1 (V, U) - \dim_k
\Hom_\Lambda (V, U) + \sum_{x \in \Delta_0} d_x' d_x''
\end{gather*}
(compare~\cite{Bongartz1994}*{2.1}).

\subsection{} \label{subsect_regular}
Let $\Lambda$ be an algebra and $\bd$ a dimension vector. It can
be easily verified that $T_M \mod_\Lambda (\bd)$ is a subspace of
$\bbZ (M, M)$ for each $M \in \mod_\Lambda (\bd)$
(compare~\cite{Voigt1977}*{3.4}). Using this observation one may
show for $M \in \mod_\Lambda (\bd)$ that if $\Lambda$ is
triangular, $\gldim \Lambda \leq 2$, and $\Ext_\Lambda^2 (M, M) =
0$, then $\dim_k \bbZ (M, M) = a_\Lambda (\bd) = \sum_{\alpha \in
\Delta_1} d_{s_\alpha} d_{t_\alpha} - \sum_{\rho \in \frakR}
d_{s_\rho} d_{t_\rho}$ and $M$ is a regular point of $\mod_\Lambda
(\bd)$ (compare~\cite{BobinskiSkowronski1999a}*{Proposition~3.2}).
Observe that if $\Lambda$ is triangular and $\gldim \Lambda \leq
2$, then $a_\Lambda (\bd) = \dim \GL (\bd) - q_\Lambda (\bd)$.

\subsection{} \label{subsect_homscheme}
Let $\Lambda$ be an algebra, $M$ a $\Lambda$-module, and $\bd$ a
dimension vector. For $d \in \bbN$ we define the subset
$\calH_d^{\bd} (M)$ of $\mod_\Lambda (\bd)$ consisting of all $N$
such that $\dim_k \Hom_\Lambda (M, N) = d$. This is a (possibly
empty) locally closed subset of $\mod_\Lambda (\bd)$. It follows
from~\cite{Zwara2002a}*{Section~3.4} (see
also~\cite{Zwara2002b}*{Section~2}) that if $N \in \calH_d^{\bd}
(M)$, then $T_N \calH_d^{\bd} (M)$ is contained in the subspace
$\bbZ_M (N, N)$ of $\bbZ (N, N)$ consisting of all $Z$ with
$\dim_k \Hom_\Lambda (M, W^Z) = 2 \cdot \dim_k \Hom_\Lambda (M,
N)$.

\subsection{} \label{subsect_extscheme}
Let $\Lambda$ be an algebra, and let $\bd'$ and $\bd''$ be
dimension vectors. For $d \in \bbN$ we define the subset
$\calE_d^{\bd', \bd''}$ of $\mod_\Lambda (\bd') \times
\mod_\Lambda (\bd'')$ consisting of all $(U, V)$ such that
$\Hom_\Lambda (V, U) = 0$ and $\dim_k \Ext_\Lambda^1 (V, U) = d$.
It follows from~\cite{Crawley-BoeveySchroer2002}*{Lemma~4.3} that
$\calE_d^{\bd', \bd''}$ is a (possibly empty) locally closed
subset of $\mod_\Lambda (\bd') \times \mod_\Lambda (\bd'')$. Let
$(U, V) \in \calE_d^{\bd', \bd''}$. Then $T_{U, V} \calE_d^{\bd',
\bd''}$ is contained in the subspace of $\bbZ (U, U) \times \bbZ
(V, V)$ consisting of all $(Z', Z'')$ such that the class of
$\xi^{Z'} \circ \xi + \xi \circ \xi^{Z''}$ is zero in
$\Ext_\Lambda^2 (V, U)$ for each $[\xi] \in \Ext_\Lambda^1 (V,
U)$~\cite{Bobinski2008}*{Proposition~3.3}. In particular, if
$\Lambda$ is triangular, $\gldim \Lambda \leq 2$, $\Ext_\Lambda^2
(U, U) = 0 = \Ext_\Lambda^2 (V, V)$, and there exists an exact
sequence $\xi : 0 \to U \to M \to V \to 0$ such that either
$\Ext_\Lambda^2 (M, U) = 0$ or $\Ext_\Lambda^2 (V, M) = 0$, then
\[
\dim_k T_{U, V} \calE_d^{\bd', \bd''} \leq a (\bd') + a (\bd'') -
\dim_k \Ext_\Lambda^2 (V, U).
\]
Indeed, our assumptions imply that the map $\bbZ (U, U) \times
\bbZ (V, V) \to \Ext_\Lambda^2 (V, U)$, $(Z', Z'') \mapsto
[\xi^{Z'} \circ \xi + \xi \circ \xi^{Z''}]$, is surjective, hence
the claim follows.

\subsection{} \label{subsect_projone}
Let $\Lambda$ an algebra and $\bd$ a dimension vector. For a full
subcategory $\calC$ of $\mod \Lambda$ by $\calC (\bd)$ we denote
$\mod_\Lambda (\bd) \cap \calC$. Let $\calP$ be the full
subcategory of $\mod \Lambda$ consisting of the modules of
projective dimension at most $1$. Then $\calP (\bd)$ is an open
subset of $\mod_\Lambda (\bd)$. Moreover, if $\Lambda$ is
triangular and $\calP (\bd) \neq \varnothing$, then $\calP (\bd)$
is irreducible. Finally, if additionally $\gldim \Lambda \leq 2$,
then $\dim \calP (\bd) = a_\Lambda (\bd)$
(see~\cite{BarotSchroer2001}*{Proposition~3.1} for all the above
claims).

\subsection{} \label{subsect_directsum}
Let $\Lambda$ be an algebra and $\bd'$, $\bd''$ dimension vectors.
If $\calU' \subseteq \mod_\Lambda (\bd')$ and $\calU'' \subseteq
\mod_\Lambda (\bd'')$, then $\calU' \oplus \calU''$ denotes the
set of all $M \in \mod_\Lambda (\bd' + \bd'')$ such that $M \simeq
M' \oplus M''$ for some $M' \in \calU'$ and $M'' \in \calU''$.
Obviously, if $\calU'$ and $\calU''$ are constructible
(irreducible), then $\calU' \oplus \calU''$ is also constructible
(irreducible, respectively). Moreover,
\begin{multline*}
\dim (\calU' \oplus \calU'') = \dim \calU' + \dim \calU'' + \dim
\GL (\bd' + \bd')
\\
- \dim \GL (\bd') - \min \{ \dim_k \Hom_\Lambda (U', U'') \mid
\text{$U' \in \calU'$, $U'' \in \calU''$} \}
\\
- \dim \GL (\bd'') - \min \{ \dim_k \Hom_\Lambda (U'', U') \mid
\text{$U' \in \calU'$, $U'' \in \calU''$} \}
\end{multline*}
(see for example~\cite{Bobinski2008a}*{Lemma~3.4}).

\makeatletter
\def\@secnumfont{\mdseries} 
\makeatother

\section{Geometric bisections} \label{sect_bisection}

In order to make possible to apply the arguments used in the proof
of Proposition~\ref{prop_regone} in other situations, we start
with developing the following technical concept. If $\Lambda$ is
an algebra, then a pair $(\calX, \calY)$ of full subcategories of
$\mod \Lambda$, which are closed under direct sums, will be called
a geometric bisection of $\mod \Lambda$ if the following
conditions are satisfied:
\begin{enumerate}

\item
$\calX \vee \calY =  \mod \Lambda$,

\item
$\pd_\Lambda X \leq 1$ for each $X \in \calX$ and $\id_\Lambda Y
\leq 1$ for each $Y \in \calY$,

\item
$\Hom_\Lambda (Y, X) = 0$ and $\Ext_\Lambda^1 (X, Y) = 0$ for all
$X \in \calX$ and $Y \in \calY$,

\item
if $\bd$ is a dimension vector, then $\calX (\bd)$ and $\calY
(\bd)$ are open subsets of $\mod_\Lambda (\bd)$.

\end{enumerate}

An example of a geometric bisection is provided by $(\calF (T),
\calT (T))$ for a directing tilting module $T$ over an algebra
$\Lambda$ (which is necessarily tilted). Indeed, the first three
conditions are just the properties of directing tilting modules
presented in~\eqref{subsect_tilted}. The remaining one follows
from the upper-semicontinuity of the functions
\[
\mod_\Lambda (\bd) \ni M \mapsto \dim_k \Hom_\Lambda (T, M),
\dim_k \Ext_\Lambda^1 (T, M) \in \bbZ
\]
\cite{Crawley-BoeveySchroer2002}*{Lemma~4.3}.

Another example of a geometric bisection can be constructed for a
concealed-canonical algebra $\Lambda$. Indeed, let $\calP$,
$\calR$ and $\calQ$ be as in~\eqref{subsect_concan}, $\calX =
\calP$, and $\calY = \calR \vee \calQ$. Then it remains to show
that the last condition of the definition is satisfied. However,
the proofs of~\cite{Bobinski2008a}*{Lemmas~3.7 and~3.8} given for
canonical algebras generalize easily to the case of concealed
canonical algebras.

Let $\Lambda$ be an algebra of global dimension at most $2$ with a
geometric bisection $(\calX, \calY)$. Observe that this means in
particular that $\Lambda$ is quasi-tilted. If $\calX (\bd) \neq
\varnothing$ for a dimension vector $\bd$, then $\ol{\calX (\bd)}
= \ol{\calP (\bd)}$, hence in particular $\ol{\calX (\bd)}$ is an
irreducible component of $\mod_\Lambda (\bd)$ of dimension
$a_\Lambda (\bd)$ (see~\eqref{subsect_projone}).

\begin{proposition} \label{prop_regone}
Let $\Lambda$ be a quasi-tilted algebra with a geometric bisection
$(\calX, \calY)$. Let $\bd$ be a dimension vector such that $\calX
(\bd) \neq \varnothing$ and for each $N \in \ol{\calX (\bd)}
\setminus \calX (\bd)$ such that $\calO (N)$ is maximal in
$\ol{\calX (\bd)} \setminus \calX (\bd)$ there exists an exact
sequence $0 \to N_1 \to M \to N_2 \to 0$ with $M \in \calX (\bd)$
and $N_1 \oplus N_2 \simeq N$. Then $\ol{\calX (\bd)}$ is regular
in codimension one.
\end{proposition}

The proof basically repeats the arguments
from~\cite{Bobinski2008}*{Proof of the main result}, however some
changes are necessary, hence we include it here for completeness.

\begin{proof}
Since $\pd_\Lambda M \leq 1$ for each $M \in \calX (\bd)$, $\calX
(\bd)$ consists of points which are regular in $\ol{\calX (\bd)}$.
Thus it suffices to show that for each irreducible component
$\calV$ of $\ol{\calX (\bd)} \setminus \calX (\bd)$ there exists
an open subset $\calU$ of $\calV$ consisting of regular points of
$\ol{\calX (\bd)}$. Fix such a component $\calV$. Our assumptions
imply that there exist dimension vectors $\bd'$ and $\bd''$ such
that $\calX (\bd') \oplus \calY (\bd'')$ contains an open subset
$\calU'$ of $\calV$. Put
\begin{gather*}
p = \min \{ \pd_\Lambda N \mid N \in \calU' \}, \qquad d' = \min
\{ \dim_k \End_\Lambda (N) \mid N \in \calU' \}, .
\\
\intertext{and} %
d'' = \min \{ \dim_k \Ext_\Lambda^1 (N, N) \mid N \in \calU' \}.
\end{gather*}
Let $\calU$ be the set of all $N \in \calU'$ such that $N$ does
not belong to an irreducible component of $\ol{\calX (\bd)}
\setminus \calX (\bd)$ different from $\calV$, $\pd_\Lambda N =
p$, $\dim_k \End_\Lambda (N) = d'$, and $\dim_k \Ext_\Lambda^1 (N,
N) = d''$. Then $\calU$ is an open subset of $\calV$. We show that
$N$ is a regular point of $\ol{\calX (\bd)}$ for each $N \in
\calU$. The claim is obvious if $p \leq 1$, thus we may assume
that $\pd_\Lambda N = 2$ for each $N \in \calU'$.

We first show that for each $N \in \calU$ there exists an exact
sequence $0 \to N' \to M \to N'' \to 0$ with $M \in \calX (\bd)$,
$N' \in \calX (\bd')$, $N'' \in \calY (\bd'')$, and $N' \oplus N''
\simeq N$. Since $\calO (N)$ is maximal in $\ol{\calX (\bd)}
\setminus \calX (\bd)$, there exists an exact sequence $\xi : 0
\to N_1 \to M_0 \to N_2 \to 0$ such that $N_1 \oplus N_2 \simeq N$
and $M_0 \in \calX (\bd)$. Choose indecomposable direct summands
$N_1'$ and $N_2'$ of $N_1$ and $N_2$, respectively, such that $p
\circ \xi \circ i : 0 \to N_1' \to M' \to N_2' \to 0$ does not
split, where $p : N_1 \to N_1'$ and $i : N_2' \to N_2$ are the
canonical projection and injection, respectively. Write $N = N_0
\oplus N_1' \oplus N_2'$ and let $M = N_0 \oplus M'$. Then $N \not
\simeq M$, $\calO (N) \subseteq \ol{\calO (M)}$, and $\calO (M)
\subseteq \ol{\calO (M_0)}$. Consequently, the maximality of
$\calO (N)$ in $\ol{\calX (\bd)} \setminus \calX (\bd)$ implies
that $M \in \calX$. This immediately implies that $N_0 \oplus N_1'
\in \calX$. Moreover, $N_2' \in \calY$, since $\pd_\Lambda N = 2$
and $N_2'$ is indecomposable, hence by adding $0 \to N_0 \to N_0
\to 0 \to 0$ to $p \circ \xi \circ i$ we obtain the desired exact
sequence.

Now we prove that $N$ in regular in $\ol{\calX (\bd)}$ for each $N
\in \calU$. Let $0 \to N' \to M \to N'' \to 0$ be an exact
sequence with $M \in \calX (\bd)$, $N' \in \calX (\bd')$, $N'' \in
\calY (\bd'')$, and $N' \oplus N'' \simeq N$. Obviously, we may
assume that $N = N' \oplus N''$, i.e.\ $N \in \mod_\Lambda (\bd')
\times \mod_\Lambda (\bd'')$, where we identify $\mod_\Lambda
(\bd') \times \mod_\Lambda (\bd'')$ with
\[
\left\{ \left(
\begin{bmatrix}
U_\alpha & 0 \\ 0 & V_\alpha
\end{bmatrix}
\right)_{\alpha \in \Delta_1} \Big| \; \text{$U \in \mod_\Lambda
(\bd')$, $V \in \mod_\Lambda (\bd'')$} \right\} \subseteq
\mod_\Lambda (\bd).
\]
Let $\calU_0$ be the intersection of $\calU$ with $\mod_\Lambda
(\bd') \times \mod_\Lambda (\bd'')$. Then $\calU_0$ is an open
subset of $(\mod_\Lambda (\bd') \times \mod_\Lambda (\bd'')) \cap
\ol{\calX (\bd)}$. Moreover, direct calculations show that
$\calU_0 \subseteq \calE_d^{\bd', \bd''}$ for $d = d'' - d' +
q_\Lambda (\bd') + q_\Lambda (\bd'') + \langle \bd', \bd''
\rangle_\Lambda$. Consequently,
\begin{multline*}
\dim_k T_N ((\mod_\Lambda (\bd') \times \mod_\Lambda (\bd'')) \cap
\ol{\calX (\bd)}) = \dim_k T_N \calU_0
\\
\leq \dim_k T_{N', N''} \calE_d^{\bd', \bd''} \leq a_\Lambda
(\bd') + a_\Lambda (\bd'') - \dim_k \Ext_\Lambda^2 (N'', N')
\end{multline*}
according to~\eqref{subsect_extscheme}.

Now we show that $T_N \ol{\calX (\bd)} \leq a (\bd)$ and this will
finish the proof. Indeed, the above inequality implies that
\begin{multline*}
T_N \ol{\calX (\bd)} \leq a_\Lambda (\bd') + a_\Lambda (\bd'') -
\dim_k \Ext_\Lambda^2 (N'', N')
\\
+ \dim_k \bbZ (N', N'') + \dim_k \bbZ (N'', N')
\end{multline*}
and direct calculations prove the claim.
\end{proof}

\makeatletter
\def\@secnumfont{\mdseries} 
\makeatother

\section{Proof of the main result} \label{sect_proof}

We start with presenting the following consequence
of~\cite{BobinskiSkowronski1999a}*{Theorems~1 and~2}. If $\Lambda$
is a tame quasi-tilted algebra and $\bd$ is the dimension vector
of an indecomposable $\Lambda$-module, then $\mod_\Lambda (\bd)$
has at most two irreducible components and each irreducible
component of $\mod_\Lambda (\bd)$ has dimension $a_\Lambda (\bd)$.
Moreover, if $\mod_\Lambda (\bd)$ is irreducible, then
$\mod_\Lambda (\bd)$ is normal, hence in particular regular in
codimension one. On the other hand, if $\mod_\Lambda (\bd)$ is not
irreducible, then we may assume that $\bd$ is sincere and one of
the following conditions is satisfied:
\begin{enumerate}

\item \label{case1}
$\bd = \bh + \bd'$ for nonzero connected dimension vectors $\bh$
and $\bd'$ with disjoint supports such that $q_\Lambda (\bh) = 0$,
$q_\Lambda (\bd') = 1$, and $d_x' \leq 1$ each vertex $x$ of the
Gabriel quiver of $\Lambda$,

\item \label{case2}
$\bd = \bh' + \bh''$ for connected dimension vectors $\bh'$ and
$\bh''$ such that $q_\Lambda (\bh') = 0 = q_\Lambda (\bh'')$,
$\langle \bh', \bh'' \rangle_\Lambda = 1$, and $\langle \bh'',
\bh' \rangle_\Lambda = 0$.

\end{enumerate}
We remark that in both cases $q_\Lambda (\bd) = 1$, hence there
exists a uniquely determined (up to isomorphism) indecomposable
$\Lambda$-module $X$ of dimension vector $\bd$. We study now the
above described cases. In the below considerations we will use
freely information about $\mod \Lambda$ and $\mod_\Lambda (\bd)$
for $\Lambda$ and $\bd$ as above obtained
in~\cites{BobinskiSkowronski1999a, BobinskiSkowronski1999b}.

First assume that the condition~\eqref{case1} is satisfied. In
this case $\Lambda$ is either tilted or concealed canonical, hence
in particular it possesses a geometric bisection $(\calX, \calY)$.
Up to duality, we may assume that $X \in \calX$. A discussion of
this case presented in~\cite{BobinskiSkowronski1999a}*{Section~5}
shows that under this assumption one of the irreducible components
of $\mod_\Lambda (\bd)$ is $\ol{\calX (\bd)}$. Moreover, if $\calO
(M)$ is maximal in $\ol{\calX (\bd)}$ for $M \in \mod_\Lambda
(\bd)$, then $M \in \calX (\bd)$. Indeed, one shows that the
maximal $\GL (\bd)$-orbits in $\ol{\calX (\bd)}$ are precisely the
maximal $\GL (\bd)$-orbits in $\mod_\Lambda (\bd)$ which are
contained in $\calX (\bd)$. The proof of this fact requires a more
detailed description of $\mod \Lambda$ as given
in~\cite{BobinskiSkowronski1999a}*{Section~5} and we leave it to
the reader as an exercise. In particular, \eqref{subsect_deg}
implies that for each $N \in \ol{\calX (\bd)} \setminus \calX
(\bd)$ such that $\calO (N)$ is maximal in $\ol{\calX (\bd)}
\setminus \calX (\bd)$ there exists an exact sequence $0 \to N_1
\to M \to N_2 \to 0$ with $M \in \calX (\bd)$ and $N_1 \oplus N_2
\simeq N$. Consequently, $\ol{\calX (\bd)}$ is regular in
codimension one according to Proposition~\ref{prop_regone}. The
other irreducible component of $\mod_\Lambda (\bd)$ in the
case~\eqref{case1} is (isomorphic to) $\mod_\Lambda (\bh) \times
\mod_\Lambda (\bd')$. The results
of~\cite{BobinskiSkowronski1999a} imply that both factors are
normal, hence regular in codimension one, and the claim follows
for this component.

Now we study the case~\eqref{case2}. The irreducible components of
$\mod_\Lambda (\bd)$ appearing in this case are described
in~\cite{BobinskiSkowronski1999b}. This description implies that
one of the irreducible components of $\mod_\Lambda (\bd)$ is
$\ol{\calO (X)}$. Moreover, $X$ is directing, hence the claim
follows from~\cite{BobinskiZwara2006}*{Theorem~1.1}. Thus it
remains to study the other irreducible component $\calV$ of
$\mod_\Lambda (\bd)$. We describe this situation more precisely.
For $p, q, r, s, t \in \bbN_+$ let $\Delta (p, q, r, s, t)$ be the
quiver
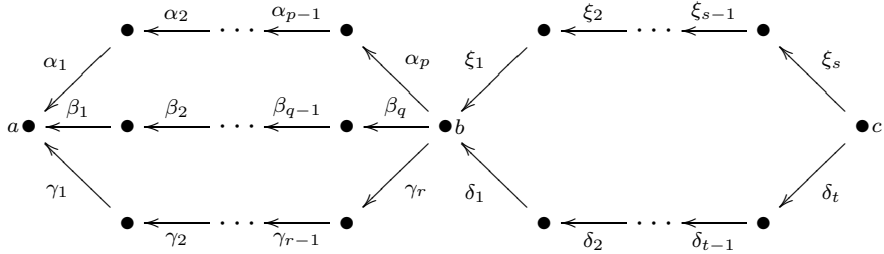
\[
\xymatrix{& \bullet \ar[ld]_{\alpha_1} & \cdots \ar[l]_-{\alpha_2}
& \bullet \ar[l]_-{\alpha_{p - 1}} & & \bullet \ar[ld]_{\xi_1} &
\cdots \ar[l]_-{\xi_2} & \bullet \ar[l]_-{\xi_{s - 1}}
\\
\bullet \save*+!R{\scriptstyle a} \restore & \bullet
\ar[l]_{\beta_1} & \cdots \ar[l]_-{\beta_2} & \bullet
\ar[l]_-{\beta_{q - 1}} & \bullet \save*+!L{\scriptstyle b}
\restore \ar[lu]_{\alpha_p} \ar[l]_{\beta_q} \ar[ld]^{\gamma_r} &
& & & \bullet \save*+!L{\scriptstyle c} \restore \ar[lu]_{\xi_s}
\ar[ld]^{\delta_t}
\\
& \bullet \ar[lu]^{\gamma_1} & \cdots \ar[l]^-{\gamma_2} & \bullet
\ar[l]^-{\gamma_{r - 1}} & & \bullet \ar[lu]^{\delta_1} & \cdots
\ar[l]^-{\delta_2} & \bullet \ar[l]^-{\delta_{t - 1}}}
\]
and
\begin{multline*}
\frakR (p, q, r, s, t) = \{ \alpha_1 \cdots \alpha_p - \beta_1
\cdots \beta_q + \gamma_1 \cdots \gamma_r,
\\
\alpha_p \xi_1, \beta_q \xi_1 \cdots \xi_s - \beta_q \delta_1
\cdots \delta_t, \gamma_r \delta_1 \}.
\end{multline*}
Then either $\Lambda$ or $\Lambda^{\op}$ is isomorphic to $k
\Delta (p, q, r, s, t) / \frakR (p, q, r, s, t)$ for some $p, q,
r, s, t \in \bbN_+$. Observe that if at least one of the numbers
$p$, $q$, $r$ equals $1$, then $(\Delta (p, q, r, s, t), \frakR
(p, q, r, s, t))$ is not a bound quiver. In this case a
corresponding bound quiver is obtained by removing an appropriate
arrow (this arrow may be not uniquely determined if at least two
of the numbers $p$, $q$, $r$ equal $1$) and modifying $\frakR (p,
q, r, s, t)$ accordingly
(see~\cite{BobinskiSkowronski2001}*{Section~2} for the list of
bound quivers obtained in this way). This operation induces an
isomorphism of the corresponding representation varieties, hence
we may work with the above family of algebras in order to reduce
the number of considered cases.

Thus assume that $\Lambda = k \Delta / \frakR (p, q, r, s, t)$ for
some $p, q, r, s, t \in \bbN_+$, where we put $\Delta = \Delta (p,
q, r, s, t)$. Let $\Delta'$ ($\Delta''$) be the minimal convex
subquiver of $\Delta$ containing $a$ and $b$ ($b$ and $c$,
respectively). If $\bh'$ and $\bh''$ are as in~\eqref{case2}, then
\[
h'_x =
\begin{cases}
1 & x \in \Delta_0'
\\
0 & x \not \in \Delta_0'
\end{cases}
\qquad \text{and} \qquad h''_x =
\begin{cases}
1 & x \in \Delta_0''
\\
0 & x \not \in \Delta_0''
\end{cases}.
\]
For $\lambda \in k \setminus \{ 0, 1 \}$ we define the
representation $H' (\lambda)$ of $\Delta$ of dimension vector
$\bh'$ by
\[
H' (\lambda)_\alpha =
\begin{cases}
\lambda & \alpha = \alpha_1,
\\
\lambda + 1 & \alpha = \beta_1,
\\
1 & \alpha \in \Delta_1', \; \alpha \neq \alpha_1, \beta_1,
\\
0 & \alpha \in \Delta_1''.
\end{cases}
\]
Similarly, we define the representation $H'' (\lambda)$ of
$\Delta$ of dimension vector $\bh''$ for each $\lambda \in k
\setminus \{ 0 \}$. Moreover, for each $i \in [1, p]$ we define
the representation $H' (\alpha_i)$ of $\Delta$ of dimension vector
$\bh'$ by
\[
H' (\alpha_i)_\alpha =
\begin{cases}
0 & \alpha = \alpha_i,
\\
1 & \alpha \in \Delta_1', \; \alpha \neq \alpha_i,
\\
0 & \alpha \in \Delta_1''.
\end{cases}
\]
Similarly, we define the representations $H' (\beta_i)$, $i \in
[1, q]$, (in this case $-1$ has to appear in the definition) and
$H' (\gamma_i)$, $i \in [1, r]$, of dimension vector $\bh'$, and
the representations $H'' (\xi_i)$, $i \in [1, s]$, and $H''
(\delta_i)$, $i \in [1, t]$, of dimension vector $\bh''$. Let
$\calU$ be the union of $\calO (H' (u) \oplus H'' (v))$, $u \in (k
\setminus \{ 0, 1 \}) \cup \Delta_1'$, $v \in (k \setminus \{ 0
\}) \cup \Delta_1''$. Then $\calV$ is the closure of $\calU$.

Now we show that $\calV = \mod_\Lambda (\bh') \oplus \mod_\Lambda
(\bh'')$. Since $\calU \subseteq \mod_\Lambda (\bh') \oplus
\mod_\Lambda (\bh'')$ and $\mod_\Lambda (\bh') \oplus \mod_\Lambda
(\bh'')$ is irreducible, it suffices to show that $\mod_\Lambda
(\bh') \oplus \mod_\Lambda (\bh'')$ is closed. However, it can be
easily checked that $M \in \mod_\Lambda (\bh') \oplus \mod_\Lambda
(\bh'')$ if and only if $\rk [(M_\alpha^{\tr})_{\substack{\alpha
\in \Delta_1 \\ s_\alpha = b}}] \leq 1$, $\rk
[(M_\alpha)_{\substack{\alpha \in \Delta_1
\\ t_\alpha = b}}] \leq 1$, and $M_{\alpha'} M_{\alpha''} = 0$ for all
$\alpha', \alpha'' \in \Delta_1$ with $s_{\alpha'} = b =
t_{\alpha''}$.

As the next step we show that $\dim (\calV \setminus \calU) \leq
\dim \calV - 2$. Let $\calU'$ be the union of $\calO (H' (u))$, $u
\in (k \setminus \{ 0, 1 \}) \cup \Delta_1'$. Then $\mod_\Lambda
(\bh') \setminus \calU'$ consists of $M \in \mod_\Lambda (\bh')$
such that $M_{\alpha_i} = M_{\beta_j} = M_{\gamma_l} = 0$ for some
$i \in [1, p]$, $j \in [1, q]$, and $l \in [1, r]$. Thus $\dim
(\mod_\Lambda (\bh') \setminus \calU') = \dim \mod_\Lambda (\bh')
- 2$. Similarly, $\dim (\mod_\Lambda (\bh'') \setminus \calU'') =
\dim \mod_\Lambda (\bh'') - 2$, where $\calU''$ is the union of
$\calO (H'' (v))$, $v \in (k \setminus \{ 0 \}) \cup \Delta_1''$.
Since
\[
\calV \setminus \calU = (\calU' \oplus \mod_\Lambda (\bh'')) \cup
(\mod_\Lambda (\bh') \oplus \calU''),
\]
the inequality $\dim (\calV \setminus \calU) \leq \dim \calV - 2$
follows from~\eqref{subsect_directsum}.

As the last step of the proof we show that $\calU$ consists of
regular points of $\calV$. Let $S$ be the simple $\Lambda$-module
at $b$, i.e.\ $S_x = \delta_{x, b} k$ for $x \in \Delta_0$, where
$\delta_{x, y}$ is the Kronecker delta. Then $\calU \subseteq
\calH_1^{\bd} (S)$ and $\calH_1^{\bd} (S)$ contains an open subset
of $\calV$. Consequently, $T_M \calV \subseteq \bbZ_S (M, M)$ for
each $M \in \calU$~\eqref{subsect_homscheme}. However, if $M = H'
(u) \oplus H'' (v)$ for $u \in (k \setminus \{ 0, 1 \}) \cup
\Delta_1'$ and $v \in (k \setminus \{ 0 \}) \cup \Delta_1''$, then
direct calculations show that
\begin{multline*}
\bbZ_S (M, M) = \bbZ (H' (u), H' (u)) \oplus \bbZ (H'' (v), H''
(v))
\\
\oplus \bbZ (H' (u), H'' (v)) \oplus \bbB (H'' (v), H' (u)).
\end{multline*}
By applying formulas from~\eqref{subsect_cocyle} we get $\dim_k
\bbZ_S (M, M) \leq a_\Lambda (\bd)$, and this finishes the proof.


\bibsection


\begin{biblist}

\bib{AssemSimsonSkowronski2006}{book}{
   author={Assem, I.},
   author={Simson, D.},
   author={Skowro{\'n}ski, A.},
   title={Elements of the Representation Theory of Associative Algebras. Vol. 1},
   series={London Math. Soc. Stud. Texts},
   volume={65},
   publisher={Cambridge Univ. Press},
   place={Cambridge},
   date={2006},
   pages={x+458},
   isbn={978-0-521-58423-4},
   isbn={978-0-521-58631-3},
   isbn={0-521-58631-3},
}

\bib{Bakke1988}{article}{
   author={Bakke, {\O}.},
   title={Some characterizations of tilted algebras},
   journal={Math. Scand.},
   volume={63},
   date={1988},
   number={1},
   pages={43--50},
   issn={0025-5521},
}

\bib{BarotSchroer2001}{article}{
   author={Barot, M.},
   author={Schr{\"o}er, J.},
   title={Module varieties over canonical algebras},
   journal={J. Algebra},
   volume={246},
   date={2001},
   number={1},
   pages={175--192},
   issn={0021-8693},
}

\bib{BenderBongartz2003}{article}{
   author={Bender, J.},
   author={Bongartz, K.},
   title={Minimal singularities in orbit closures of matrix pencils},
   journal={Linear Algebra Appl.},
   volume={365},
   date={2003},
   pages={13--24},
   issn={0024-3795},
}

\bib{Bobinski2002}{article}{
   author={Bobi{\'n}ski, G.},
   title={Geometry of decomposable directing modules over tame algebras},
   journal={J. Math. Soc. Japan},
   volume={54},
   date={2002},
   number={3},
   pages={609--620},
   issn={0025-5645},
}

\bib{Bobinski2008a}{article}{
   author={Bobi{\'n}ski, G.},
   title={Geometry of regular modules over canonical algebras},
   journal={Trans. Amer. Math. Soc.},
   volume={360},
   date={2008},
   number={2},
   pages={717--742},
   issn={0002-9947},
}

\bib{Bobinski2008}{article}{
   author={Bobi{\'n}ski, G.},
   title={Orbit closures of directing modules are regular in codimension one},
   eprint={arXiv:0712.1246},
   status={preprint}
}

\bib{BobinskiSkowronski1999a}{article}{
   author={Bobi{\'n}ski, G.},
   author={Skowro{\'n}ski, A.},
   title={Geometry of modules over tame quasi-tilted algebras},
   journal={Colloq. Math.},
   volume={79},
   date={1999},
   number={1},
   pages={85--118},
   issn={0010-1354},
}

\bib{BobinskiSkowronski1999b}{article}{
   author={Bobi{\'n}ski, G.},
   author={Skowro{\'n}ski, A.},
   title={Geometry of directing modules over tame algebras},
   journal={J. Algebra},
   volume={215},
   date={1999},
   number={2},
   pages={603--643},
   issn={0021-8693},
}

\bib{BobinskiSkowronski2001}{article}{
   author={Bobi{\'n}ski, G.},
   author={Skowro{\'n}ski, A.},
   title={Selfinjective algebras of Euclidean type with almost regular nonperiodic Auslander-Reiten components},
   journal={Colloq. Math.},
   volume={88},
   date={2001},
   number={1},
   pages={93--120},
   issn={0010-1354},
}

\bib{BobinskiSkowronski2002}{article}{
   author={Bobi{\'n}ski, G.},
   author={Skowro{\'n}ski, A.},
   title={Geometry of periodic modules over tame concealed and tubular algebras},
   journal={Algebr. Represent. Theory},
   volume={5},
   date={2002},
   number={2},
   pages={187--200},
   issn={1386-923X},
}

\bib{BobinskiZwara2002}{article}{
   author={Bobi{\'n}ski, G.},
   author={Zwara, G.},
   title={Schubert varieties and representations of Dynkin quivers},
   journal={Colloq. Math.},
   volume={94},
   date={2002},
   number={2},
   pages={285--309},
   issn={0010-1354},
}

\bib{BobinskiZwara2006}{article}{
   author={Bobi{\'n}ski, G.},
   author={Zwara, G.},
   title={Normality of orbit closures for directing modules over tame algebras},
   journal={J. Algebra},
   volume={298},
   date={2006},
   number={1},
   pages={120--133},
   issn={0021-8693},
}

\bib{Bongartz1983}{article}{
   author={Bongartz, K.},
   title={Algebras and quadratic forms},
   journal={J. London Math. Soc. (2)},
   volume={28},
   date={1983},
   number={3},
   pages={461--469},
   issn={0024-6107},
}

\bib{Bongartz1991}{article}{
   author={Bongartz, K.},
   title={A geometric version of the Morita equivalence},
   journal={J. Algebra},
   volume={139},
   date={1991},
   number={1},
   pages={159--171},
   issn={0021-8693},
}

\bib{Bongartz1994}{article}{
   author={Bongartz, K.},
   title={Minimal singularities for representations of Dynkin quivers},
   journal={Comment. Math. Helv.},
   volume={69},
   date={1994},
   number={4},
   pages={575--611},
   issn={0010-2571},
}

\bib{Bongartz1996}{article}{
   author={Bongartz, K.},
   title={On degenerations and extensions of finite-dimensional modules},
   journal={Adv. Math.},
   volume={121},
   date={1996},
   number={2},
   pages={245--287},
   issn={0001-8708},
}

\bib{Crawley-Boevey1988}{article}{
   author={Crawley-Boevey, W. W.},
   title={On tame algebras and bocses},
   journal={Proc. London Math. Soc. (3)},
   volume={56},
   date={1988},
   number={3},
   pages={451--483},
   issn={0024-6115},
}

\bib{Crawley-BoeveySchroer2002}{article}{
   author={Crawley-Boevey, W. W.},
   author={Schr{\"o}er, J.},
   title={Irreducible components of varieties of modules},
   journal={J. Reine Angew. Math.},
   volume={553},
   date={2002},
   pages={201--220},
   issn={0075-4102},
}

\bib{Drozd1980}{collection.article}{
   author={Drozd, Yu. A.},
   title={Tame and wild matrix problems},
   book={
      title={Representation Theory. II},
      series={Lecture Notes in Math.},
      volume={832},
      editor={Dlab, V.},
      editor={Gabriel, P.},
      publisher={Springer},
      place={Berlin},
   },
   date={1980},
   pages={242--258},
}

\bib{Eisenbud1995}{book}{
   author={Eisenbud, D.},
   title={Commutative Algebra with a View toward Algebraic Geometry},
   series={Grad. Texts in Math.},
   volume={150},
   publisher={Springer},
   place={New York},
   date={1995},
   pages={xvi+785},
   isbn={0-387-94268-8},
   isbn={0-387-94269-6},
}

\bib{Geiss1996}{collection.article}{
   author={Gei{\ss}, Ch.},
   title={Geometric methods in representation theory of finite-dimensional algebras},
   book={
      title={Representation Theory of Algebras and Related Topics},
      series={CMS Conf. Proc.},
      volume={19},
      editor={Bautista, R.},
      editor={Mart{\'{\i}}nez-Villa, R.},
      editor={de la Pe{\~n}a, J. A.},
      publisher={Amer. Math. Soc.},
      place={Providence},
   },
   date={1996},
   pages={53--63},
}

\bib{GeissSchroer2003}{article}{
   author={Geiss, Ch.},
   author={Schr{\"o}er, J.},
   title={Varieties of modules over tubular algebras},
   journal={Colloq. Math.},
   volume={95},
   date={2003},
   number={2},
   pages={163--183},
   issn={0010-1354},
}

\bib{GeissSchroer2005}{article}{
   author={Geiss, Ch.},
   author={Schr{\"o}er, J.},
   title={Extension-orthogonal components of preprojective varieties},
   journal={Trans. Amer. Math. Soc.},
   volume={357},
   date={2005},
   number={5},
   pages={1953--1962},
   issn={0002-9947},
}

\bib{HappelReiten1999}{article}{
   author={Happel, D.},
   author={Reiten, I.},
   title={Hereditary categories with tilting object},
   journal={Math. Z.},
   volume={232},
   date={1999},
   number={3},
   pages={559--588},
   issn={0025-5874},
}

\bib{HappelReitenSmalo1996}{article}{
   author={Happel, D.},
   author={Reiten, I.},
   author={Smal{\o}, S. O.},
   title={Tilting in abelian categories and quasitilted algebras},
   journal={Mem. Amer. Math. Soc.},
   volume={120},
   date={1996},
   number={575},
   pages={viii+ 88},
   issn={0065-9266},
}

\bib{HappelRingel1982}{article}{
   author={Happel, D.},
   author={Ringel, C. M.},
   title={Tilted algebras},
   journal={Trans. Amer. Math. Soc.},
   volume={274},
   date={1982},
   number={2},
   pages={399--443},
   issn={0002-9947},
}

\bib{Kerner1989}{article}{
   author={Kerner, O.},
   title={Tilting wild algebras},
   journal={J. London Math. Soc. (2)},
   volume={39},
   date={1989},
   number={1},
   pages={29--47},
   issn={0024-6107},
}

\bib{Kraft1984}{book}{
   author={Kraft, H.},
   title={Geometrische Methoden in der Invariantentheorie},
   series={Aspects Math.},
   volume={D1},
   publisher={Vieweg},
   place={Braunschweig},
   date={1984},
   pages={x+308},
   isbn={3-528-08525-8},
}

\bib{KraftRiedtmann1986}{collection.article}{
   author={Kraft, H.},
   author={Riedtmann, Ch.},
   title={Geometry of representations of quivers},
   book={
      title={Representations of Algebras},
      series={London Math. Soc. Lecture Note Ser.},
      volume={116},
      editor={Webb, P.},
      publisher={Cambridge Univ. Press},
      place={Cambridge},
   },
   date={1986},
   pages={109--145},
}

\bib{Kunz1985}{book}{
   author={Kunz, E.},
   title={Introduction to Commutative Algebra and Algebraic Geometry},
   publisher={Birkh\"auser Boston},
   place={Boston},
   date={1985},
   pages={xi+238},
   isbn={3-7643-3065-1},
}

\bib{LenzingMeltzer1996}{article}{
   author={Lenzing, H.},
   author={Meltzer, H.},
   title={Tilting sheaves and concealed-canonical algebras},
   book={
      title={Representation Theory of Algebras},
      series={CMS Conf. Proc.},
      volume={18},
      editor={Bautista, R.},
      editor={Mart{\'{\i}}nez-Villa, R.},
      editor={de la Pe{\~n}a, J. A.},
      publisher={Amer. Math. Soc.},
      place={Providence},
   },
   date={1996},
   pages={455--473},
}

\bib{LenzingdelaPena1997}{article}{
   author={Lenzing, H.},
   author={de la Pe{\~n}a, J. A.},
   title={Wild canonical algebras},
   journal={Math. Z.},
   volume={224},
   date={1997},
   number={3},
   pages={403--425},
   issn={0025-5874},
}

\bib{LenzingdelaPena1999}{article}{
   author={Lenzing, H.},
   author={de la Pe{\~n}a, J. A.},
   title={Concealed-canonical algebras and separating tubular families},
   journal={Proc. London Math. Soc. (3)},
   volume={78},
   date={1999},
   number={3},
   pages={513--540},
   issn={0024-6115},
}

\bib{LenzingSkowronski1996}{article}{
   author={Lenzing, H.},
   author={Skowro{\'n}ski, A.},
   title={Quasi-tilted algebras of canonical type},
   journal={Colloq. Math.},
   volume={71},
   date={1996},
   number={2},
   pages={161--181},
   issn={0010-1354},
}

\bib{delaPena1993}{article}{
   author={de la Pe{\~n}a, J. A.},
   title={Tame algebras with sincere directing modules},
   journal={J. Algebra},
   volume={161},
   date={1993},
   number={1},
   pages={171--185},
   issn={0021-8693},
}

\bib{Richmond2001}{article}{
   author={Richmond, N. J.},
   title={A stratification for varieties of modules},
   journal={Bull. London Math. Soc.},
   volume={33},
   date={2001},
   number={5},
   pages={565--577},
   issn={0024-6093},
}

\bib{Ringel1984}{book}{
   author={Ringel, C. M.},
   title={Tame Algebras and Integral Quadratic Forms},
   series={Lecture Notes in Math.},
   volume={1099},
   publisher={Springer},
   place={Berlin},
   date={1984},
   pages={xiii+376},
   isbn={3-540-13905-2},
}

\bib{Schroer2004}{article}{
   author={Schr{\"o}er, J.},
   title={Varieties of pairs of nilpotent matrices annihilating each other},
   journal={Comment. Math. Helv.},
   volume={79},
   date={2004},
   number={2},
   pages={396--426},
   issn={0010-2571},
}

\bib{Skowronski1996}{article}{
   author={Skowro{\'n}ski, A.},
   title={On omnipresent tubular families of modules},
   book={
      title={Representation Theory of Algebras},
      series={CMS Conf. Proc.},
      volume={18},
      editor={Bautista, R.},
      editor={Mart{\'{\i}}nez-Villa, R.},
      editor={de la Pe{\~n}a, J. A.},
      publisher={Amer. Math. Soc.},
      place={Providence},
   },
   date={1996},
   pages={641--657},
}

\bib{Skowronski1998}{article}{
   author={Skowro{\'n}ski, A.},
   title={Tame quasi-tilted algebras},
   journal={J. Algebra},
   volume={203},
   date={1998},
   number={2},
   pages={470--490},
   issn={0021-8693},
}

\bib{SkowronskiZwara2003}{article}{
   author={Skowro{\'n}ski, A.},
   author={Zwara, G.},
   title={Derived equivalences of selfinjective algebras preserve singularities},
   journal={Manuscripta Math.},
   volume={112},
   date={2003},
   number={2},
   pages={221--230},
   issn={0025-2611},
}

\bib{Voigt1977}{book}{
   author={Voigt, D.},
   title={Induzierte Darstellungen in der Theorie der endlichen, algebraischen Gruppen},
   series={Lecture Notes in Math.},
   volume={592},
   publisher={Springer},
   place={Berlin},
   date={1977},
   pages={iv+413},
   isbn={3-540-08251-4},
}

\bib{Zwara2000}{article}{
   author={Zwara, G.},
   title={Degenerations of finite-dimensional modules are given by extensions},
   journal={Compositio Math.},
   volume={121},
   date={2000},
   number={2},
   pages={205--218},
   issn={0010-437X},
}

\bib{Zwara2002b}{article}{
   author={Zwara, G.},
   title={Unibranch orbit closures in module varieties},
   journal={Ann. Sci. \'Ecole Norm. Sup. (4)},
   volume={35},
   date={2002},
   number={6},
   pages={877--895},
   issn={0012-9593},
}

\bib{Zwara2002a}{article}{
   author={Zwara, G.},
   title={Smooth morphisms of module schemes},
   journal={Proc. London Math. Soc. (3)},
   volume={84},
   date={2002},
   number={3},
   pages={539--558},
   issn={0024-6115},
}

\bib{Zwara2003}{article}{
   author={Zwara, G.},
   title={An orbit closure for a representation of the Kronecker quiver with bad singularities},
   journal={Colloq. Math.},
   volume={97},
   date={2003},
   number={1},
   pages={81--86},
   issn={0010-1354},
}

\bib{Zwara2005b}{article}{
   author={Zwara, G.},
   title={Orbit closures for representations of Dynkin quivers are regular in codimension two},
   journal={J. Math. Soc. Japan},
   volume={57},
   date={2005},
   number={3},
   pages={859--880},
   issn={0025-5645},
}

\bib{Zwara2005}{article}{
   author={Zwara, G.},
   title={Regularity in codimension one of orbit closures in module varieties},
   journal={J. Algebra},
   volume={283},
   date={2005},
   number={2},
   pages={821--848},
   issn={0021-8693},
}

\bib{Zwara2006}{article}{
   author={Zwara, G.},
   title={Singularities of orbit closures in module varieties and cones over rational normal curves},
   journal={J. London Math. Soc. (2)},
   volume={74},
   date={2006},
   number={3},
   pages={623--638},
   issn={0024-6107},
}

\bib{Zwara2007}{article}{
   author={Zwara, G.},
   title={Codimension two singularities for representations of extended Dynkin quivers},
   journal={Manuscripta Math.},
   volume={123},
   date={2007},
   number={3},
   pages={237--249},
   issn={0025-2611},
}

\end{biblist}
\end{document}